\documentclass[a4paper,12pt]{article}

\usepackage[utf8]{inputenc}
\usepackage{amsmath}
\usepackage{amsfonts}
\usepackage{amssymb}
\usepackage{amsthm}
\usepackage{caption}
\usepackage[english]{babel}
\usepackage{fontenc}
\usepackage{graphicx}
\usepackage{xcolor}
\usepackage{hyperref}

\usepackage{algorithm}
 \usepackage{algorithmic}

\def\R {{\mathbb R}}
\def\E {{\mathbb E}}

\newcommand{\X}{X_t}

\newcommand{\de}{\partial}

\newcommand{\uu}{\mathrm{u}}

\title{\bf Fokker-Planck analysis of superresolution microscopy images}
\author{Mario Annunziato 
\thanks{Dipartimento di Fisica ``E.~R.~Caianiello'', Universit\`a degli Studi di Salerno, 
  Via G. Paolo II 132, 84084 Fisciano, Italy. e-mail: mannunzi@unisa.it, mannunzi@am-research.it}
\and Alfio Borz{\`i}
\thanks{Institut f\"ur Mathematik, Universit\"at W\"urzburg, Emil-Fischer-Strasse 30,
97074 W\"urzburg, Germany. e-mail: alfio.borzi@mathematik.uni-wuerzburg.de}}
\date{\today}

\begin{document}

 \maketitle

 \begin{abstract}
 A method for the analysis of superresolution microscopy images is presented. 
 This method is based on the analysis of stochastic trajectories of particles 
 moving on the membrane of a cell with the assumption that this 
 motion is determined by the properties of this membrane. Thus, the purpose 
 of this method is to recover the structural properties of the membrane 
 by solving an inverse problem governed by the Fokker-Planck
equation related to the stochastic trajectories. Results of numerical experiments 
demonstrate the ability of the proposed method to reconstruct the 
potential of a cell membrane by using synthetic data similar those captured by superresolution
microscopy of luminescent activated proteins.
 \end{abstract}
 

 \section{Introduction}
 \label{sec:intro}

The pioneering works \cite{Betzig2006,Hell2005,Moerner2003} mark 
the development of the revolutionary superresolution microscopy (SRM) 
that allows to go beyond the Abbe limit for conventional light microscopy \cite{Abbe1873}. 
The SRM method consists in labelling the molecules moving 
on a biological support with fluorophores and then in sampling 
the microscopic images of the activated fluorescent molecules.

Observation of the frames of the sampled SRM microscopic images have suggested that the motion of the molecules could be modelled by a stochastics Langevin
equation \cite{Saffman1975,Schuss2009}. Specifically, it appears that an adequate model of the observed trajectories of 2-dimensional images is given by the following stochastic differential equation (SDE) \cite{Holcman2015}: 
\begin{align}
d\X &=  b(\X)  \, dt +   \sigma(\X) \, dW_t  
\label{eq:sde2}\\
X_{t_0} & = X_0,  \label{eq:ic2}
\end{align}
where $b$ represents the drift, $\sigma$ is the dispersion coefficient, and 
$\X\in\R^2$ denotes the position of the observed molecule at time $t$. 
In this framework, it is well-known that the drift and dispersion coefficients satisfy
$$
\lim_{t \to s} \E\left[ \,  \frac{\X - X_s}{t-s} \, | \,  X_s=z  \right] = b(z) , \qquad
\lim_{t \to s} \E\left[ \, \frac{ | \X - X_s |^2 }{t-s} \,  | \,  X_s=z  \right] = \sigma^2(z) ,
$$
where the expected values are computed with respect to the process 
having value $z$ at $t=s$; the operator $\E[\cdot \, | \,  X_s=z ]$ denotes 
averaging w.r.t the measure of the trajectories conditioned to be at $z$ at time $s$.

The formulas above suggest that suitable approximations to $b$ 
and $\sigma$ can be obtained by tracking single molecules; 
see, e.g., \cite{IMAP2015,Holcman2015}. However, this 
approach may suffer of the highly fluctuating values 
of the trajectories and the difficulty of discerning between 
different molecules that come closer than the resolution limit. 

For this reason, already in \cite{MarioAlfioJSC} the authors have 
pursued an alternative strategy that allows to build a robust methodology for estimation of the drift based on the observation of 
ensemble of trajectories. Our approach is built upon the assumption 
that this ensemble is driven by a velocity field (the drift), 
given by a potential velocity field $U(x)$, as follows: 
\begin{equation}\label{eq:drift}
 b(x;U)=- \nabla U(x).
\end{equation}
Moreover, one assumes 
a constant diffusion coefficient whose value is chosen consistently with 
estimates of laboratory measurement \cite{diffusion}.

Further, in agreement with our statistical approach based on ensembles, 
we focus on the evolution of the probability density function (PDF) of the positions of the molecules (not on the single trajectories) whose evolution is governed by the Fokker-Planck (FP) problem given by \cite{Risken1996,Schuss2009}: 
 \begin{align}
 & \partial_t f(x,t) - \nabla \cdot \left( \nabla U(x) \, f (x,t)\right) - \frac{\sigma^2 }{2} \, \Delta f (x,t)=0 , \quad (x,t) \in Q \label{FPn1} \\
 & F (f) \cdot \hat n = 0, \qquad (x,t) \in \Sigma \label{eq:zero_flux} , \\
 & f(x,0)  = f_0(x) , \qquad x \in \Omega \label{eq:initialFPMFaux},
 \end{align}
 where $Q=\Omega\times (0,T) $ and $\Sigma =\partial \Omega\times (0,T) $. 
 In this formulation, $f(x,t)$ represents the PDF of a particle at $x \in \R^2$ at time $t$, $\nabla U(x)$ 
is the Cartesian gradient of the potential $U$, $f_0$ is the initial 
density, and $\Delta$ is the two-dimensional Laplace operator. Notice 
that we require zero-flux boundary conditions, where $F(f)(x,t)$ is the 
following flux of probability
 \begin{equation}
 F(f)(x,t)=\frac{\sigma^2 }{2} \, \nabla f(x,t)  - b(x;U) \,  f(x,t) .
 \label{Flux1}
 \end{equation}
We choose zero-flux boundary conditions since they reasonably 
model the situation where a similar number of particles enters and exits the domain; 
see, e.g., \cite{MarioAlfioJCAM,Schuss2009}.

Our proposal is to construct a FP-based imaging modality that is based on the 
formulation of an inverse problem for $U$ and the observation of a time sequence, in a time interval $[0,T]$, of numerical PDFs (two-dimensional histograms), 
which are obtained from a uniform binning of SRM particles' positions. 
We denote this input data with $f_d(x,t)$ which is a piecewise constant function.
In this setting the initial condition is given by $f_0(x)=f_d(x,0)$. 

This proposal is similar to that in our previous work \cite{MarioAlfioJSC}. 
However, in \cite{MarioAlfioJSC} the assumption of interacting 
particles was made that resulted in a nonlinear FP model and very 
involved and CPU time demanding calculations.

At the continuous level, our FP-based imaging tool is formulated 
as the following inverse problem: 
\begin{align}\label{OPCn1}
\min J(f,U)&:=\frac{1}{2} \int_\Omega \int_0^T (f(x,t) - f_d(x,t))^2\; dx\,dt \notag \\
& + \frac{\xi}{2} \int_\Omega  (f(x,T) - f_d(x,T))^2\; dx \notag \\
& + \frac{\alpha}{2} \int_\Omega (  |U (x)|^2 + | \nabla U (x)|^2 )\; dx , \\
\textit{s.t.} \quad \partial_t f (x,t) & + \nabla \cdot 
\left[ b(x,;U) f(x,t)   \right]  - \frac{\sigma^2}{2} \, \Delta f(x,t)  = 0 , \quad 
\textrm{ in } Q \notag \\
 f(x,0)  &= f_0(x) \textrm{ in }  \Omega , \qquad F(f)\cdot \hat{n} =0 \textrm{ on } 
 \Sigma, \notag
\end{align}
with the given initial- and boundary conditions for the FP equation, 
and $\alpha, \xi >0$.

In this problem, the objective functional $J$ is defined as the weighted sum 
of a space-time bestfit term 
$\int_\Omega \int_0^T (f(x,t) - f_d(x,t))^2\; dx\,dt$, and 
at final time $\int_\Omega  (f(x,T) - f_d(x,T))^2\; dx$, and of a suitable `energy' of 
the potential $ \|U \|^2 = \int_\Omega (|U (x)|^2 + | \nabla U (x)|^2 )\; dx$, 
which corresponds to the square of the $H^1(\Omega)$ norm of $U$. Notice that this formulation allows to avoid any differentiation of the data and 
makes possible to choose the binning size and, in general, the measurement setting, 
independently of any choice of parameters that are required in the numerical 
solution of the optimization problem.

Our second main concern in determining the potential $U$ is to provide a 
measure of uncertainty, and thus of reliability, of its reconstruction. Statistically,
this is achieved by many repetition of the same experiment, that could not be
feasible for (short) living cells.
However, inspired by the so-called model predictive control (MPC) scheme \cite{GruenePannek}
already used for optimal control problems \cite{reviewMarioAlfio,MarioAlfioJCAM}, 
we propose a novel procedure 
to quantify the uncertainty of the estimation of $U$ 
by using the data of a single experiment.

Our methodology is to consider 
a sequence of non-parametric inverse problems 
like \eqref{OPCn1} defined on time 
windows $(t_k,t_{k+1})$, $k=0,\ldots,K-1$, that represent an uniform partition in $K$ 
subintervals of the time interval $[0,T]$. Therefore a statistical analysis 
can be performed on the set of the corresponding $K$ solutions for $U$ 
that are obtained in the subintervals.

For development and validation, we consider images of cell's membrane structures 
(actin, cytoskeleton), as expression of potentials, that is pixel grey values, with which 
we generate our synthetic data. In particular, we use an image 
of actin from a cytoskeleton obtained with a Platinium-replica electron microscopy
\cite{dent_et_al,leb_et_al}.

\begin{figure}[ht!]
 \centering
  \includegraphics[scale=0.8]{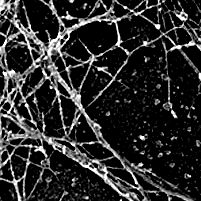}\\

\hspace*{-5mm}
\includegraphics[scale=0.356]{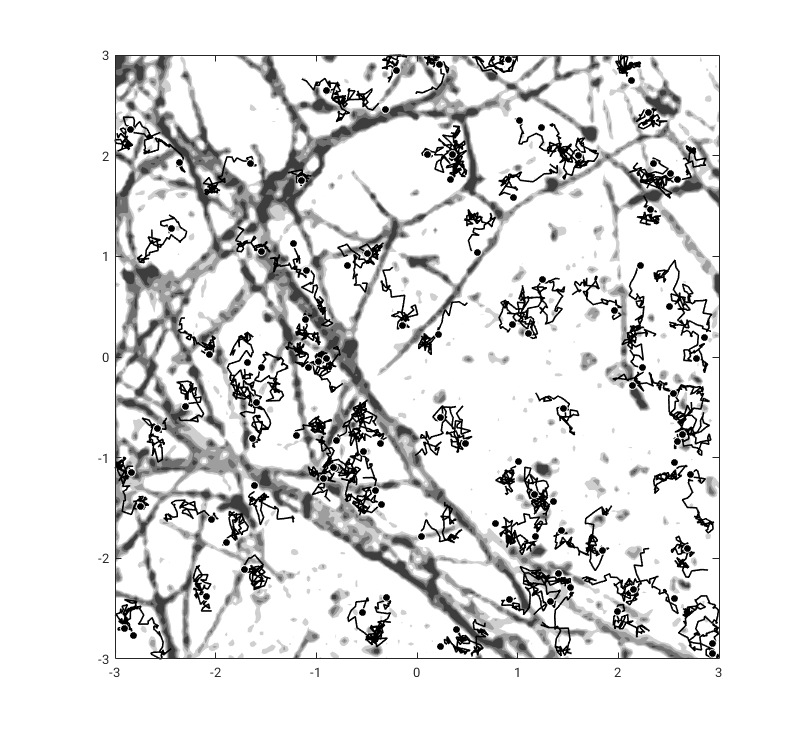}

 \caption{A picture of actin from a cytoskeleton as cell membrane potential (close up)  (Courtesy of Koch Institute \cite{Koch}) (left); a few simulated trajectories of particles (black dots)
 on the membrane (in reverse colors).}
 \label{fig:protein_U}
\end{figure}

With this images taken as gray-level representation of potential functions, we perform Monte Carlo simulation of motion of particles to generate images of molecules at different time instants, thus constructing the datasets representing the output of measurements. This setting is illustrated in Figure \ref{fig:protein_U}, where the image of actin \cite{Koch} and a plot of a 
few trajectories of the corresponding stochastic motion of the particles in this 
potential are shown.

Once the synthetic measurement data is constructed, we perform a pre-processing step on this data to construct the 
numerical PDF required in our method and solve our inverse problem to find the estimated-measured potential $U$. 
The latter is compared with that one used in the MC simulation, by a measure of similarity
based on the pixel cross-correlation between the two images.

 In Section \ref{sec:ident_probl}, 
 we discuss a numerical methodology for solving our FP-based reconstruction method for the potential $U$. 
  In Section \ref{sec: expdesign}, we provide all details of our experimental setting and introduce some analysis tools for 
 determining the accuracy of the proposed reconstruction. 
 In Section \ref{sec:numExp}, we validate our reconstruction method, 
 and use our uncertainty quantification procedure. 
   In Section \ref{sec:analysisRes}, we investigate the resolution of the proposed FP image reconstruction as optical instrument.  
 A section of conclusion and acknowledgements completes this work.

 \section{Numerical methodology}
 \label{sec:ident_probl}

\begin{figure}[!ht]
\includegraphics[scale=0.45]{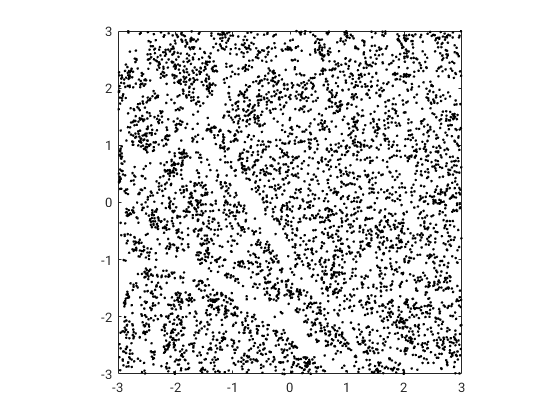}
\includegraphics[scale=0.4]{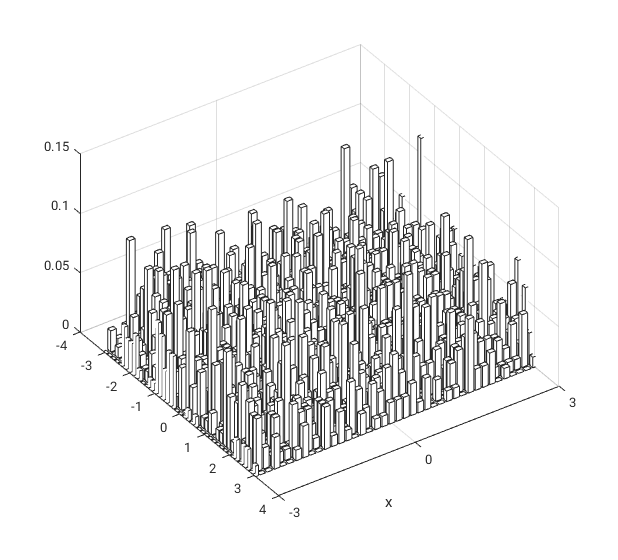}

  \caption{A frame of particles (left) and the corresponding histogram $f_d(x,t)$ on a mesh of 40x40 bins for a fixed time, from simulated data.}
 \label{fig1frames}
\end{figure}

Our aim is to reconstruct the potential $U$ from the data consisting 
of a temporally sampled SRM images of the positions of particles 
subject to this potential; see Figure \ref{fig1frames} (left), for 
a schematic snap-shot of this data. This image is subject to a 
pre-processing binning procedure in order to construct histograms 
by counting the number of particles in a regular square
partition of $\Omega$. The height of an histogram is proportional to 
the number of particles in a bin of the domain. 
This procedure for the image at time $t$ defines 
the histogram function $f_d(\cdot, t)$; see Figure \ref{fig1frames} (right).

In order to illustrate our numerical framework, we introduce 
the potential-to-state map $U \mapsto f=S(U)$, that is, the map 
that associates to a given $U \in H^1(\Omega)$ the unique solution 
to our FP problem \eqref{FPn1}-\eqref{eq:initialFPMFaux},
with given initial condition $f_0$. For the analysis of well-posedness 
and regularity of the map $S$ we refer to \cite{Koerner2022}.

Next, we remark that with the map $S$, we can define 
the reduced objective functional $\hat J(U):=J(S(U),U)$ and consider the 
equivalent formulation of \eqref{OPCn1} given by 
\begin{equation}
\min_{ U \in H^1(\Omega)} \hat J(U),
\label{OPCreduced}
\end{equation}
which has the structure of an unconstrained optimization problem. 
Thanks to the regularity of $S$ and the quadratic structure 
of $J$, existence of an optimal $U$ can be stated by well known 
techniques; see, e.g., \cite{Lions1971}.

Further, since $S$ and $J$ are Fr\'echet differentiable, 
it is possible to characterize an optimal $U$ as the 
solution to the following first-order optimality condition 
$$
\nabla_U \hat J(U) =0, 
$$
where $\nabla_U  \hat J(U)$ denotes the so-called reduced gradient  \cite{BorziSchulz2012}.

In the Lagrange framework, this condition results in 
the following optimality system: 
\begin{align}\label{OptSys1}
\partial_t f (x,t) & + \nabla \cdot 
\left[ b(x;U) f(x,t) \right]  - \frac{\sigma^2}{2} \, \Delta f(x,t)  = 0 , \notag \\
 f(x,0)  &= f_0(x) \textrm{ in }  \Omega , \qquad F(f)\cdot \hat{n} =0 \textrm{ on } 
 \partial \Omega\times (0,T] , \notag \\
 \partial_t p(x,t) & + \frac{\sigma^2}{2} \Delta p(x,t) + 
\nabla p(x,t) \cdot b(x;U) = f(x,t)-f_d(x,t) , \\
p(x,T)& =- \xi \, (f(x,T) - f_d(x,T)) \textrm{ in }   \Omega , \qquad \de_{\hat n} p(x,t) =0 \textrm{ on }  \partial \Omega \times (0,T] , \notag \\
\alpha \, U(x) & - \alpha \, \Delta U(x) - 
 \int_0^T \nabla \cdot \left(  f(x,t) \, \nabla p(x,t) \right) dt=0     \textrm{ in }   \Omega  ,   \notag \\
 & \qquad \partial_{\hat n} U=0 ,  \qquad \textrm{ on } \partial \Omega , \notag 
 \end{align}
where $p$ denotes the adjoint variable, which is governed by a 
backward adjoint FP equation. 

The last equation in \eqref{OptSys1} is the so-called the optimality condition equation, and the Neumann 
boundary condition $\partial_{\hat n} U=0$ is our modelling choice. One 
can show that its left-hand side represents the $L^2$ gradient along the FP 
differential constraint with respect to $U$ of the objective functional. We have 
\begin{equation}\label{eq:gradL2}
 \nabla_U \hat J(U)(x) := \alpha \, U(x) - \alpha \, \Delta U(x) - 
 \int_0^T \nabla \cdot \left(  f(x,t) \, \nabla p(x,t) \right) dt .
\end{equation}

Our approach for solving our FP optimization problem 
\eqref{OPCn1} is based on the nonlinear conjugate gradient (NCG) method; 
see, e.g., \cite{BorziSchulz2012}. This is an iterative method that resembles 
the standard CG scheme and requires to estimate the 
reduced gradient $\nabla_U \hat J(U)$ at each iteration. 

In order to illustrate the NCG method, we start with a discussion 
on the construction of the gradient. For a given $U^n$ obtained after $n$ iterations, 
we solve the FP equation and its adjoint, and use 
\eqref{eq:gradL2} to assemble the $L^2$ gradient. However, since the 
potential is sought in $H^1(\Omega)$, we need to obtain the $H^1$ gradient 
that satisfies the following relation 
\begin{equation}
\left( \nabla_U \hat J(U)\vert_{H^1},\delta U \right)_{H^1} = \left( \nabla_U \hat J(U)\vert_{L^2},\delta U \right)_{L^2} ,
\end{equation}
where $(\cdot, \cdot)$ denotes the $L^2(\Omega)$ scalar product. 

Now, using the definition of the $H^1$ inner product, we obtain
\begin{equation}
\int_\Omega \! \left[ \nabla_U \hat J(U)\vert_{H^1} \cdot \delta U (x) + \nabla_x \nabla_U \hat J(U)\vert_{H^1} \cdot  \nabla_x \delta U(x) \right] \,dx 
= \int_\Omega \! \nabla_U \hat J(U)\vert_{L^2} \, \delta U (x) \, dx,
\end{equation}
which must hold for all the test functions $\delta U \in H^1(\Omega)$. Therefore 
we obtain 
\begin{equation}\label{eq:gradH1}
-\Delta \left[ \nabla_U \hat J(U)\vert_{H^1} \right] + \left[ \nabla_U \hat J(U)\vert_{H^1} \right] = \nabla_U \hat  J(U)\vert_{L^2},
\end{equation}
with the boundary conditions $\frac{\de }{\de \hat n} \nabla_U \hat J(U)\vert_{H^1} = 0$ on $\de\Omega$; see \cite{MarioAlfioJSC} for more details. In Algorithm \ref{algo:CalculateGradient} our procedure for computing the gradient is given.

\begin{algorithm}
\caption{Calculate $H^1$ gradient.}
\begin{algorithmic}
\REQUIRE control $U(x)$, $f_0(x)$, $f_d(x,t)$. 
\ENSURE reduced gradient $\nabla_U \hat J (U)\vert_{H^1}$
\STATE Solve forward the FP equation with inputs: $f_0(x)$, $U(x)$
\STATE Solve backward the adoint FP equation with inputs: $U(x)$, $f(x,t)$
\STATE Assemble the $L^2$ gradient $\nabla_U \hat J(U)\vert_{L^2}$ using 
\eqref{eq:gradL2}. 
\STATE Compute the $H^1$ gradient $\nabla_U \hat J(U)\vert_{H^1}$ 
solving \eqref{eq:gradH1}. 
\RETURN $\nabla_U \hat J(U)\vert_{H^1}(x)$
\end{algorithmic}
\label{algo:CalculateGradient}
\end{algorithm}

In this algorithm, the FP problem and its optimization 
FP adjoint are approximated by the exponential Chang-Cooper scheme 
and the implicit BDF2 method; see \cite{MarioAlfioJCAM}.

Now, we can discuss the NCG method. The NCG iterative procedure 
is initialized with $U^0(x)=0$. We denote 
the optimization directions with ${d}^n$. In the first update, 
we have ${d}^0 = - \nabla_U \hat J(U^0)\vert_{H^1}$ and perform 
the optimization step
$$
U^{1} = U^0 + \alpha_0 \, {d}^0, 
$$
where $\alpha_0$ is obtained by a backtracking linesearch procedure. 
After the first step, in the NCG method the descent direction is defined 
as a linear combination of the new gradient and the past direction as follows: 
$$
{d}^{n} = -  \nabla_U \hat J (U^n)\vert_{H^1} + \beta_{n-1} \, {d}^{n-1}, 
$$
where $\beta_{-1} =0$, and $\beta_{n-1} =\| \nabla_U \hat J(U^n)\vert_{H^1} \|^2 /(d^{n-1}\cdot ( \nabla_U \hat J(U^{n})\vert_{H^1} - \nabla_U \hat J(U^{n-1})\vert_{H^1}) )$, that is, 
the Dai-Yuan formula; see, e.g., \cite{BorziSchulz2012}. 

The tolerance $tol$ and the maximum number of iterations 
$n_{\max}$ are used for termination criteria. Summarizing, in Algorithm \ref{algo:OptAlgo} we present the NCG procedure. \\
For our numerical experiments these algorithms have been implemented
with object oriented programming in C++, by using the numerical
libraries Armadillo \cite{arma}, Lapack \cite{lapa} and SuperLU \cite{slu}.

\begin{algorithm}
\caption{Nonlinear conjugate gradient (NCG) method}
\begin{algorithmic}
\REQUIRE $U^0(x)\equiv 0$, $f_0(x)$, $f_d(x,t)$
\ENSURE Optimal control $U(x)$ and corresponding state $f(x,t)$
\STATE $n = 0$
\STATE Assemble gradient ${g}^0=\nabla_U \hat J(U^0)\vert_{H^1} $ using Algorithm \ref{algo:CalculateGradient}; set ${d}^0=-{g}^0$.
\WHILE{$\|{g}^n\|_{H^1}>tol$ \AND $n<n_{\max}$}
\STATE Use linesearch to determine $\alpha_n$
\STATE Update control: $U^{n+1} = U^n + \alpha_n \, {d}^n$
\STATE Compute the gradient ${g}^{n+1}=\nabla_U \hat J(U^{n+1})\vert_{H^1} $ using Algorithm \ref{algo:CalculateGradient}
\STATE Calculate the new descent direction ${d}^{n+1} = - {g}^{n+1} + \beta_{n} \, {d}^{n}$
\STATE Set $n = n+1$
\ENDWHILE
\RETURN $U^n(x)$
\end{algorithmic}
\label{algo:OptAlgo}
\end{algorithm}

\section{Experimental design and analysis tools}
\label{sec: expdesign}
In a real application, the input data $f_d$ is given by frames 
of a recorded sequence of a SRM experiment, and the desired output is the reconstructed potential denoted with $U_r$. In our case, we construct 
this data based on a sample potential $U_s$ that determines the drift 
function in our stochastic model \eqref{eq:sde2}. Thus we generate our 
frames of synthetic data first by time-integrating this SDE in the chosen interval $[0,T]$, 
and choosing the initial positions of the particles randomly uniformly distributed. Next, the positions of the particles at different times are collected
in a sequence 2-dimensional bins that results in the sequence of 
distributions $f_d(x,t_\ell)$, $\ell = 1, \ldots, L$, where $L$ 
is the length of the resulting time sequence of frames. With this 
preparation, we apply our algorithm to obtain $U_r$, which represents 
the proposed reconstruction of $U_s$. A comparison between these 
two potentials allows to validate the accuracy of our reconstruction 
method (see below). 

We choose a domain $\Omega=[-3,3]\times[-3,3] $, and $U_s$ 
corresponds to Figure \ref{fig:protein_U} (left), where 
the values of $U_s$ in $\Omega$ correspond 
to the gray scale pixel values of the picture mapped in $[0,1]$. 
With this $U_s$, we perform a stochastic simulation of $N_p$ particles for a time horizon $T$, and diffusion amplitude $\sigma$. 
The particles trajectories given by \eqref{eq:sde2} 
with \eqref{eq:drift}, are computed with the Euler-Maruyama scheme with 
a time step $\tau=10^{-3}$, which results in a number of $L$ frames. In this 
simulation, reflecting barriers for the stochastic motion are implemented. 
We remark that for the following  calculations we are going to consider a relatively small value of density of particles; 
see \cite{Lukes,diffusion}.

Next, we perform a binning of the positions of the particles at each frame 
to construct $f_d$. Hence, we consider a uniform partition of $\Omega$ with non-overlapping squares; see, Figure \ref{fig1frames} for a 
plot of particles in $\Omega$ at a given time and the corresponding $f_d$. 
Notice that $f_d$ is irregular, nevertheless we do not perform smoothing of this data. The sequence of $f_d$ values enter in our bestfit functional in \eqref{OPCn1}. 

Once we have computed $U_r$ with our optimization 
procedure, we aim at providing a quantification of its uncertainty. Thus, 
we compute the following normalized cross-correlation factor 
between the reconstructed potential $U_r$ and the one used to 
generate the synthetic data $U_s$. We have 
\begin{equation}
 \label{eq:cross-corr}
 cc(U_r,U_s)=\dfrac{U_r\cdot U_s}{|U_r| \, |U_s|}.
\end{equation}
In this formula, $U_r$ and $U_s$ are considered as vectors and $\cdot$ 
represents the scalar product. Therefore if $cc=1$ we have 
that the two potentials match perfectly, whereas if its value is close to $0$ the 
two potentials are dissimilar, and $0.5$ it is poor. Notice that cross-correlation is commonly used in medical imaging and biology; see, e.g., \cite{lebbink_et_al,rigort_at_al,Zao_et_al}.

Clearly, one could consider many repetition of the simulation of the motion of the particles with the same initial condition and make the final binning 
on the average of the resulting frames. This procedure 
would result in a less fluctuating $f_d(x,t_\ell)$ that allows a better
reconstruction. However, this scenario 
seems difficult to realize in the real laboratory setting of a living cell. 
On the other hand, in SRM imaging is possible to visualize the motion of the 
particles on a cell membrane for a relatively long time ($T \gg 1$ in 
our setting) and our approach exploits this possibility considering a 
subdivision of the time interval in a number $K$ of time windows, and solving our optimization 
problem in each of these windows almost independently. This approach 
allows to improve the reconstruction $U_r$ and makes possible 
to quantify the uncertainty of the reconstruction.  
  
Now, to illustrate our approach, consider a uniform partition of $[0,T]$ in 
time windows of size $\Delta t = T/K$ with $K$ a 
positive integer. Let $t_k=k\Delta t$, $k=0,1,\ldots,K$, denote the 
starting- and end-points of the windows. At time $t_0$, 
we have the initial PDF $f_0$, and we solve our optimization problem \eqref{OPCn1}
in the interval $[t_0,t_1]$. This means that the final time is $t_1$ and 
the terminal condition for the adjoint variable is given by 
$p(x,t_1) =- \xi \, (f(x,t_1) - f_d(x,t_1))$.  The resulting potential is denoted 
with $U_1$. Thus, the solution obtained in this window provides also the PDF at $t=t_1$.

Clearly, we can repeat this procedure in the interval $(t_1,t_2)$ 
with the computed PDF at $t=t_1$ as the initial condition and 
$t_2$ as final time, to compute $U_2$. This procedure is recursive and can be repeated 
for $k=1,\ldots,K$, thus obtaining $U_k$, $k=1,\ldots,K$. 

Notice that small values of $K$ in relation to $L$ produce a rough estimate of the average potential and its standard deviation 
due to statistical fluctuations of the Monte Carlo experiments. On the other hand,
for greater values of $K$, the number of frames for each window of our approach is 
reduced when $L$ is kept fixed, thus resulting in a worsening of the reconstruction procedure.

For the purpose of our analysis, we apply a scaling of these potentials so that their point-wise values are 
in the interval $[0,1]$. This scaling is performed as follows: 
\begin{equation}
\hat U = \dfrac{U - \min(U)}{\max(U)-\min(U)}.
 \label{eq:scaling}
\end{equation}
Thereafter, we the reconstructed potential by pixel-wise average
of the $U_k$ is given by 
\begin{equation}
\langle U_r \rangle = \frac{1}{K} \sum_{k=1}^{K} \hat U_k . 
\label{Uavergae}
\end{equation}
Moreover, we can also compute the following pixel-wise standard deviation
 \begin{equation}
  \label{eq:stdev}
  sd(U_r)=\sqrt{\sum_{k=1}^{K} \dfrac{(\hat U_k - \langle U_r \rangle)^2}{K-1}}.
 \end{equation}

Next, we provide conversion formulas for 
our parameters in order to accommodate data from real laboratory experiments. 
We introduce a unit of length $\uu$ such that 
the side length $l$ of our square domain $\Omega$ 
is given by $l = 6\, \uu$, and the unit of the noise amplitude $\sigma$ is 
given by $\sqrt{\uu}/s$. In real biological experiments, the typical measure of the length $\tilde l$ of a cell membrane is given in $\mu m$. Further, the particle's diffusion constant $D=\sigma^2/2$ is given in $\mu m^2/s$, hence we have the 
correspondence $\sigma = l/\tilde l \sqrt{2 D}$ in unit $\sqrt{\uu}/s$, whose value is used 
for MC simulations.

The depth of the potential $\tilde U$ is expressed 
in unit of $K_B \bar T$, where $K_B$ is the Boltzmann constant and $\bar{T}$ the absolute temperature.
In experimental papers, the equation \eqref{eq:drift} is written with the diffusion constant
$D$ and $K_B \bar{T}$, i.e. $D \tilde U /(K_B \bar{T})$. As above, we obtain the relationship between the 
values of a potential $U$ and the scaled $\tilde U$, as $\tilde U = U (\tilde l/ l)^2/D$ in the unit of $K_B \bar{T}$.

As an illustration of the setting above, we see that in an experiment, the superresolution of an acquired image frame can reach the value of $0.02\, \mu m/pixel$. 
With an image of $500 \times 500$ pixels, we have $\tilde l = 10\; \mu m$.
The average diffusion coefficient of particles (protein molecules) observed in 
SRM imaging is estimated with $D=0.1\, \mu m^2/s$ \cite{diffusion}.
By superresolution techniques, it is possible 
to activate a density of $0.5 \div 2 / \mu m^2$ visible particles, 
which in terms of image pixels corresponds to $0.5 \div 2$ particles in a square 
of 50 pixels of side. Each frame is usually sampled at time intervals of $\delta t = 30\; ms$. 

In order to set up a consistent MC simulation of a real experiment, 
by mapping an image of a square of side 10 $\mu m$ on our domain $\Omega$, 
we get from the above mentioned formula: 
$\sigma = 6/10\sqrt{0.2} \approx 0.268\, \uu/\sqrt{s}$.

\section{Numerical validation}
\label{sec:numExp}
In this section, we discuss results of experiments 
in a setting that is close to real laboratory experiments involving 
SRM imaging. The results of these experiments demonstrate the ability of our methodology to reconstruct the potential
 from the simulations of the SRM measurements of the motion of particles on a cell's membrane.

\begin{figure}[!ht]
\centering
\includegraphics[scale=0.7]{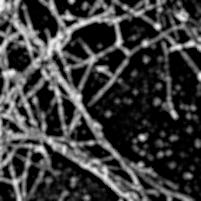}
\includegraphics[scale=1.4]{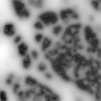}
 \caption{Left: Portion of the cytoskeleton (Courtesy of \cite{Koch}). 
 Right: reconstructed potential with the MPC scheme and $K=5$.}
 \label{fig:cito_schel}
\end{figure}

We consider a potential that corresponds 
to a portion of cytoskeleton as depicted in Figure \ref{fig:cito_schel}, with $200 \times 200$ pixels. We assume that the pixel is $50\, nm$, which corresponds to an area 
of $100\; \mu m^2$. In the figure, the white regions represent the structure of the cytoskeleton; the black ones are the 
`valleys' where the proteins are supposed to be attracted.

For the MC simulations for generating the synthetic data, we choose 
$\sigma=0.268\, \uu/\sqrt{s}$. 
This value of $\sigma$ corresponds to a diffusion constant of $D\simeq 0.1\, \mu m^2/s$.
We consider $N_p=1000$ particles, i.e. an average density of 
10 particles per $\mu m^2$. In this case, we consider a sequence of $L=3000$ frames and
$T=90$, obtained by the numerical integration of the stochastic differential equation with an integration step $\tau=30\cdot 10^{-3}\,s$. The frames have $\delta t=30\; ms$, similar to a real experiment. 
The resulting (single run) particle trajectories are collected in a binning 
process based on a mesh $\Omega$ of $50 \times 50$ bins.

For our reconstruction method, we choose a numerical 
partition of $\Omega$ of $100 \times 100$ subdivisions, corresponding to 
a mesh size of 100 $nm$. The time integration step coincides with that of the frames. For the tracking functional, we set $\alpha = 10^{-4}$ and $\xi=1$. 
Further, in the FP setting, we have $\sigma=0.7\, \uu/\sqrt{s}$. 
Notice that $\sigma$ in the FP model is chosen 
larger than the one used in the MC simulations. This choice is dictated by 
numerical convenience and it appears that it does not affect the quality of the reconstruction. The calculations are performed according to the MPC procedure 
with $K=5$ time windows. 

With this setting, we obtain the reconstructed potential shown in Figure  \ref{fig:cito_schel} (right). We see that the reconstruction is less sharp as we expected considering 
the much finer structure of the cytoskeleton and the small number of particles involved. 

\begin{figure}[!ht]
 \centering
\includegraphics[scale=0.1]{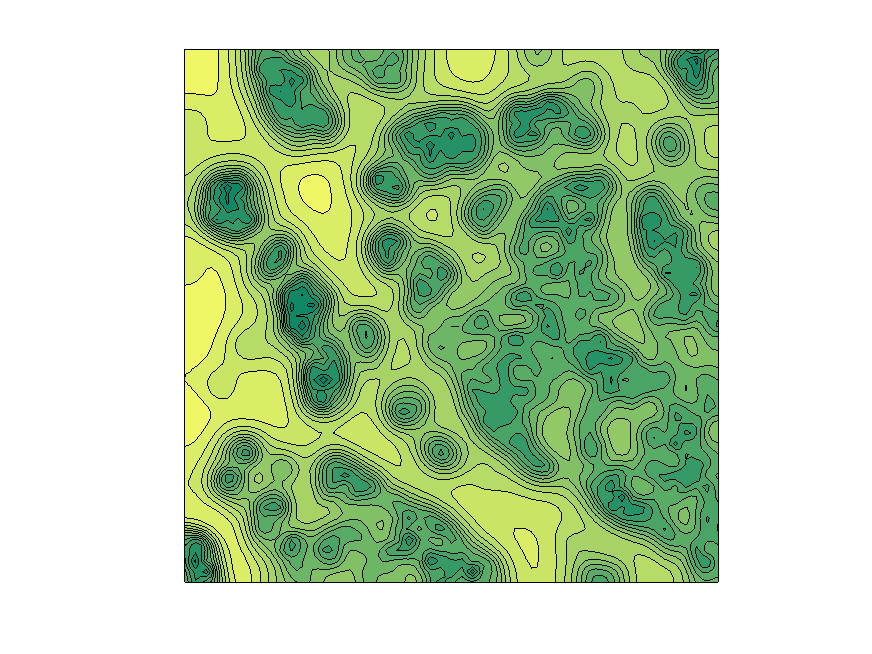}
 \hspace*{-0.75cm}
 \includegraphics[scale=0.1]{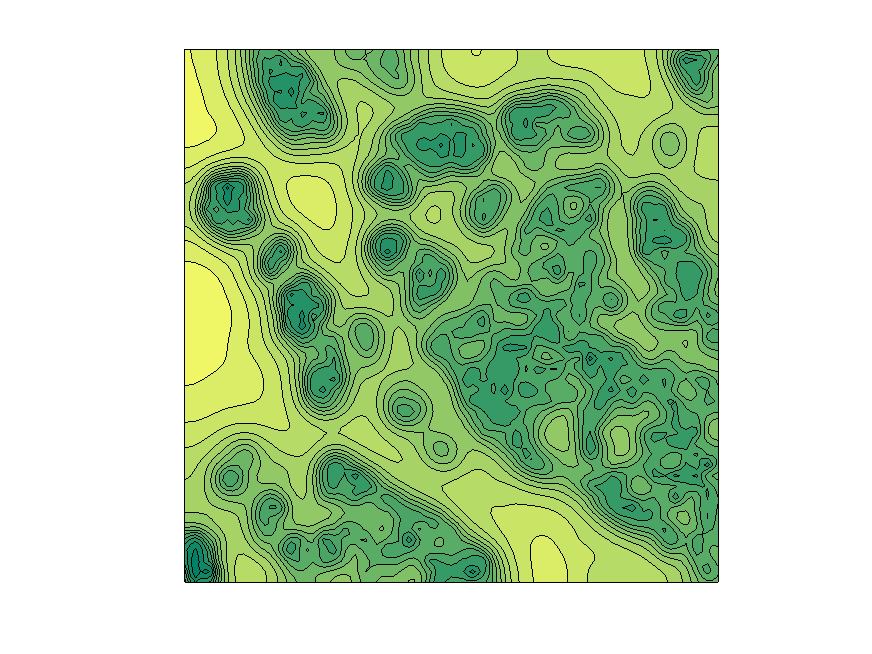}
\hspace*{-0.75cm}
\includegraphics[scale=0.1]{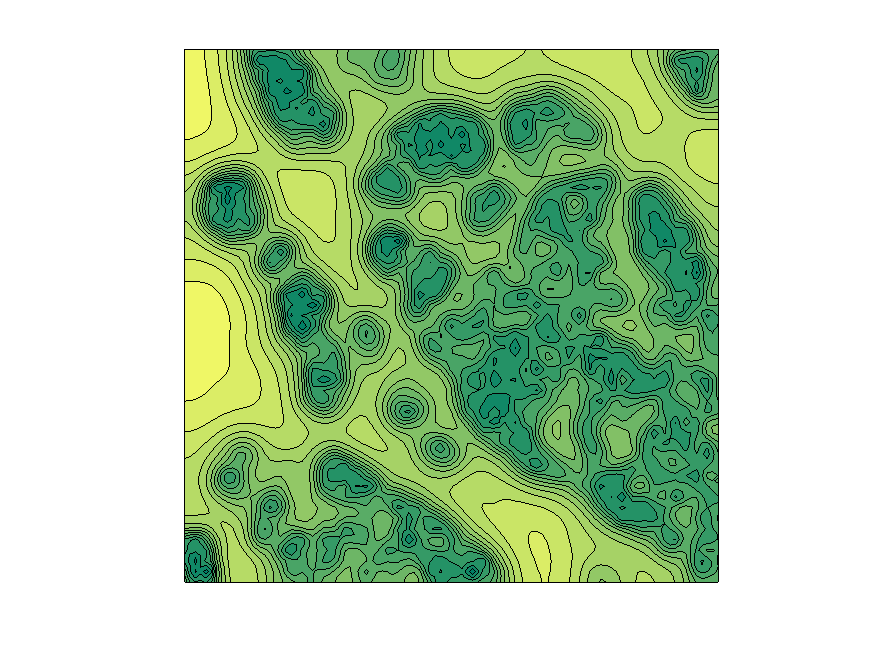}
\hspace*{-0.75cm}
\includegraphics[scale=0.1]{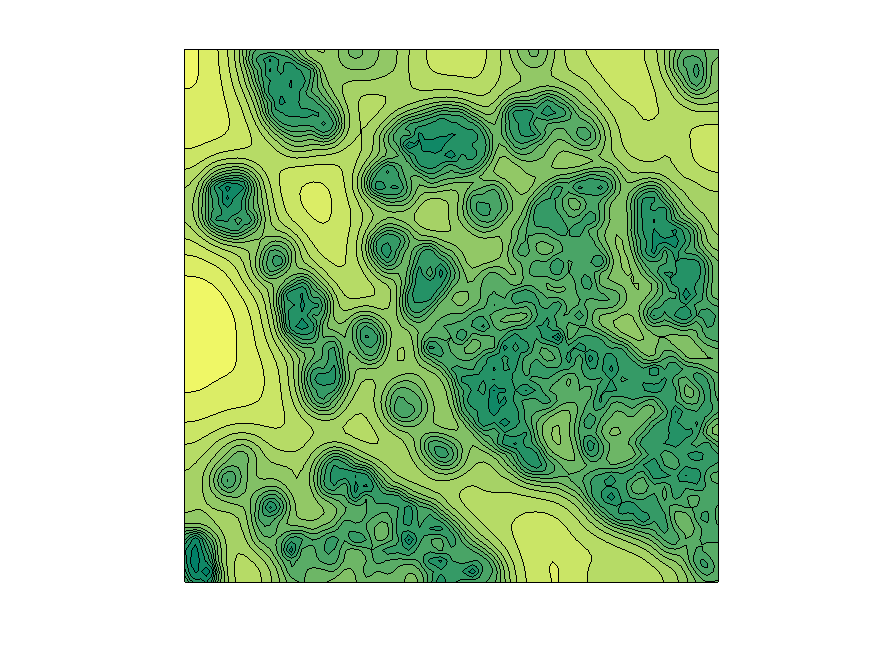}
\hspace*{-0.75cm}
\includegraphics[scale=0.1]{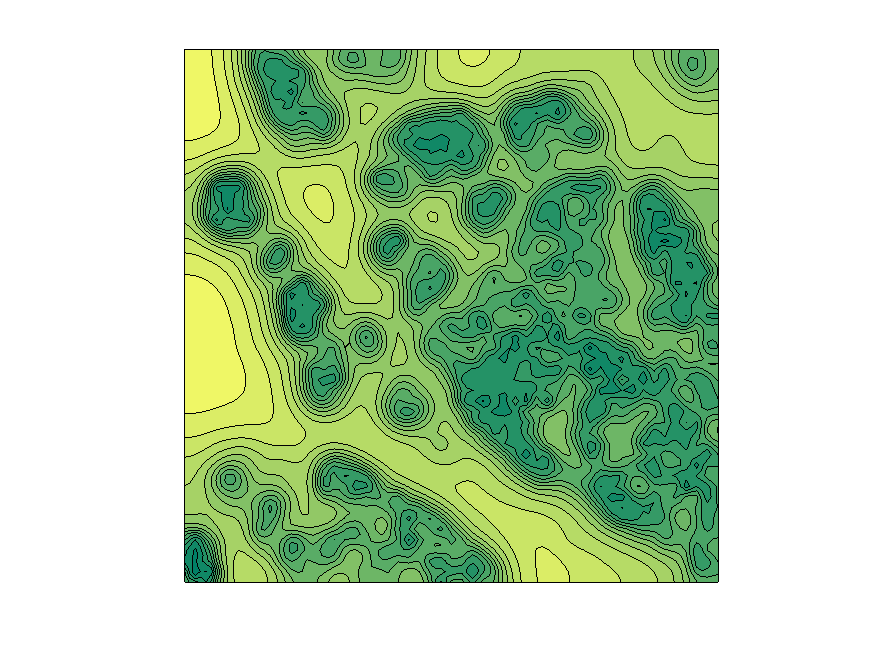}

 \caption{Sequence of the 5 (from top-left, top-right, etc.) calculated potentials obtained with the MPC procedure. 
 Cross-correlation values: 0.82, 0.81, 0.82, 0.81, 0.82}
 \label{fig:citoschel_MPC}
\end{figure}

\begin{figure}[ht!]
 \centering
\includegraphics[scale=0.2]{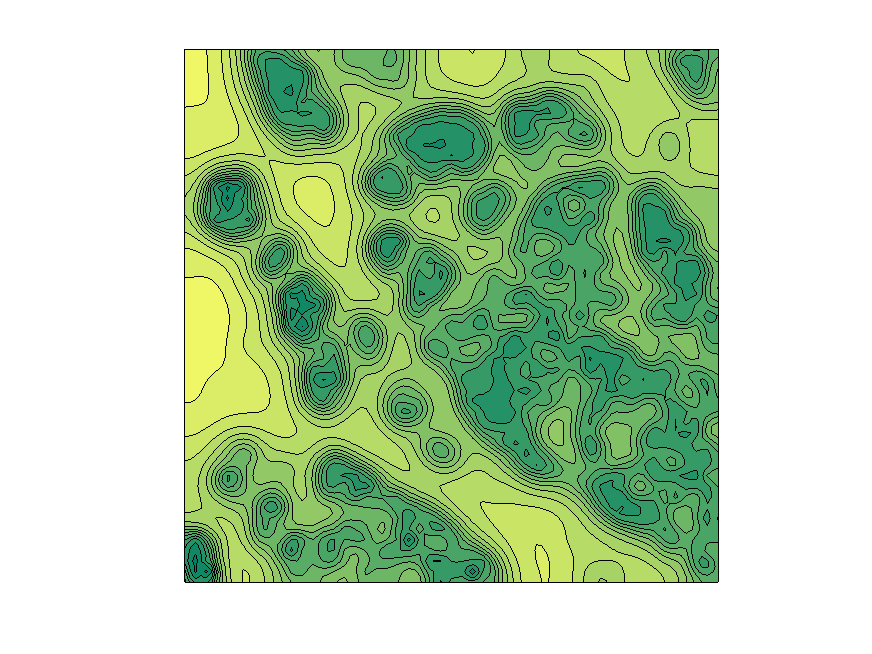}
\includegraphics[scale=0.2]{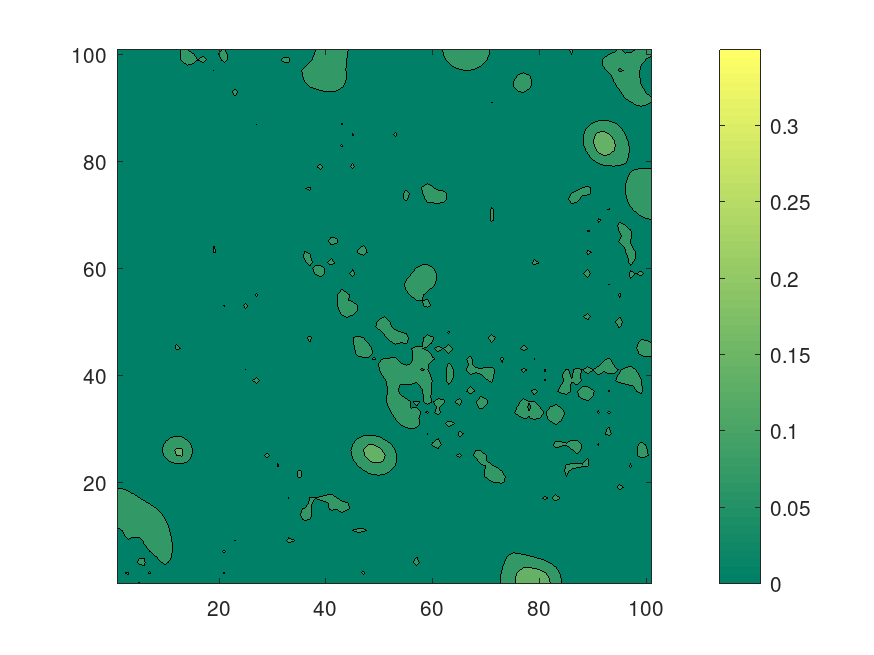}
 \caption{Left: reconstructed mean potential. Its cross-correlation value with respect to the real image is 0.82. 
 Right: standard deviation of the reconstructed potential, in level set representation.}
 \label{fig:U_citoschel}
\end{figure}

Further, in Figure \ref{fig:citoschel_MPC}, we depict the potentials obtained 
on each time window of the MPC procedure and the values of the corresponding $cc$. With these results, we have 
obtained the reconstructed potential $U_r$ in Figure  \ref{fig:cito_schel}, which we 
re-plot in level-set format in Figure \ref{fig:U_citoschel} for 
comparison. In Figure \ref{fig:U_citoschel}, we also depict the standard deviation 
that suggests that we have obtained a reliable reconstruction with 
small uncertainty. 

 \section{Resolution of FP-based image reconstruction}
 \label{sec:analysisRes}

 In this section, we investigate the optical resolution of our reconstruction method, 
 that is, try to determine a confidence value related to the scale
 at which our method can resolve variations of the potential. As a guideline, we remark that single molecule localization microscopy (SMLM) can distinguish distances of molecules of approximately $20 \, nm$ resolution. Therefore we assume this resolution range of the fluorescently labelled particles images, and attempt to quantify the smallest scale at which geometric features of the reconstructed potential $U$ can be distinguished.

 For our purpose, we consider the following `target potential', appearing as an
 alternating sequence of black and white circles (likewise those in test targets
 used for the resolution measurement of optical instruments), to synthetically generate the 
 motion data of particles. We have 
\begin{equation}
  U(x,y)=A \, \left( 1 +\cos\left(\dfrac{2 \pi}{d l}(x^2+y^2)\right) \right), \qquad 
 (x,y) \in \Omega ,
\label{testUres}
 \end{equation}
where $A$ denotes the semi-amplitude of the variation between the minimum and the
 maximum of the  potential, $l$ is the length of the side of the domain, $d$ is the 
 distance  between two peaks of the potential as a fraction of $l$. 
 
\begin{figure}[ht!]
 \centering
 \includegraphics[scale=0.35]{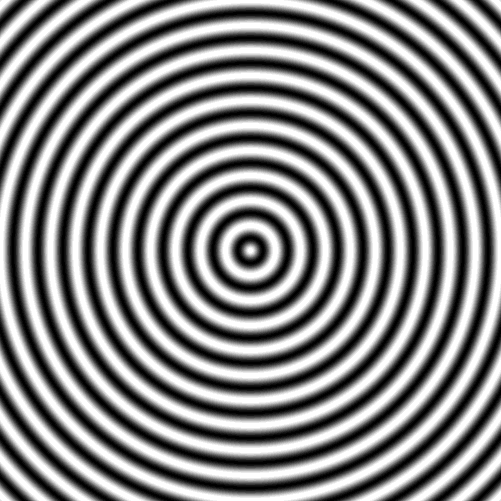}
 \includegraphics[scale=1.74]{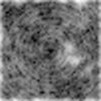}
 \caption{Left: the potential \eqref{testUres} with $A=0.05$ and $d=1/20$ (the gray scale levels spans from $U =$ 0 to 0.1). 
 Right: result of the reconstruction with the gray levels expanded to the min/max of $\langle U_r\rangle$. The cross correlation between the two images is 0.82.}
 \label{fig:calibra}
\end{figure}

 Now, we consider a single MC simulation of $500$ particles with the setting: 
 $\sigma=0.5$, $T=90$ and $L=3000$ frames, integrated with the time step $\tau = 0.03$.
 In Figure \ref{fig:calibra}, we show (left) the given potential with $A=0.05$ and $d=1/20$,
 with a gray-scale value representation conveniently adjusted for illustration pourpose.
 According to the above working hypothesis, we suppose that the pixel's width of the image is $20\, nm$
 In Figure \ref{fig:calibra} (left), we depict $U$ in a square of side of $500$ pixels,
 corresponding to $\tilde l = 10\, \mu m$. Hence, the distance between two peaks is
 $\lambda = 10 / 20 = 500\; nm$. Further, the particle's density is 
 $5$ particles per $\mu m^2$,  the diffusion coefficient $D\simeq 0.3472\, \mu m^2/s$, 
 and the potential depth, i.e. the difference between the maximum and the minimum, is $\tilde U = 0.8 \; K_B T$. 
For the reconstruction process, we use a  binning of $50 \times 50$, $\alpha=10^{-4}, \xi=1$. In the 
 numerical setting, we use a grid of $100 \times 100$ points, and $K=5$.
  Also in Figure \ref{fig:calibra} (right), we show the reconstructed potential $\langle U_r\rangle$ and notice 
 its high accuracy that is also confirmed by the high value $0.82$ of the cross correlation.
 Notice, that the quality of the reconstruction can be further improved by using 
 post-processing techniques of images.
 
Now, with the aim to define a criteria to establish the resolution measure 
for the potential,
we introduce a confidence level for the quality of the reconstructed potential
by setting a threshold for the calculated cross-correlation. 
This approach has been adopted in \cite{lebbink_et_al} for the detection of 
cellular objects from images acquired from electron tomography. 
For that purpose the authors used the threshold value of $0.5$, 
whereas in our case, we set a more strict threshold-$cc$ level equal to $0.8$. 
With this threshold, we can state 
that the test pattern depicted in Figure \ref{fig:calibra} (left) is satisfactorily 
reconstructed and determine that the resolution measure 
associated to our `imaging instrument' is $500\; nm$.
Notice that this value is affected by the value of the potential $U$ and the diffusion
$D$, and it can be further improved by changing the other 
parameters of the experiment, such as $K$ or the time $T$ of the motion sampling.

 \section{Conclusion}

 A novel method for the analysis of superresolution microscopy images was presented, 
 and applied to the reconstruction of the structure of a cell membrane potential based on observation of the motion of particles on the membrane. 

The working principle of  this method is the modeling with the Fokker-Planck equation of the ensemble 
of the stochastic trajectories of particles moving on the membrane of a cell, and the solution of an optimization problem governed by this equation, where the purpose of the optimization is to find a potential such that a least-squares bestfit term of the computed and observed particles' 
density and a Tikhonov regularization term are minimized. 

Results of numerical experiments were presented that 
demonstrated the ability of the proposed method to reconstruct the 
potential of a cell membrane by using data of superresolution
microscopy of luminescent activated proteins.

\end{document}